\documentclass[11pt]{article}
\usepackage[a4paper, margin=25mm]{geometry}
\usepackage{url,hyperref,graphicx,xcolor,amsmath,amssymb,amsthm}

\newtheorem{lemma}{Lemma}
\newtheorem{cor}{Corollary}
\newtheorem{thm}{Theorem}

\newcommand\R{\mathbb{R}}
\newcommand\E{\mathbb{E}}

\newcommand\eps{\epsilon}
\renewcommand\P{\mathbb{P}}

\DeclareMathOperator{\D}{D}
\DeclareMathOperator{\Var}{Var}
\title{An easy proof of transience of one-dimensional Rademacher walk with power weights}
\author{Stanislav Volkov}
\begin{document}
\maketitle
\begin{abstract}
 Using a version of a local limit theorem, we provide an easy proof that the Rademacher random walk with step sizes $n^\alpha$, $\alpha>1/2$ is transient.
\end{abstract}

Let $\eps_1, \eps_2, \dots$ be  a sequence of independent Rademacher random variables, i.e., $\P(\eps_i=1)=\P(\eps_i=-1)=1/2$. Let $a_1,a_2,\dots$ be a deterministic sequence of positive real numbers. Following~\cite{bhattacharya2023recurrence} and~\cite[Chapter~8]{englander2025coin}
we define the associated \emph{Rademacher random walk} as
$$ 
S_n = \sum_{k=1}^n \eps_k a_k.
$$
We are interested in recurrence vis-à-vis transience of this walk.
Since $a_n$ may in general take non-integer values, we have to modify the definition of recurrence: for a $C>0$ we say that the walk is {\em $C$-recurrent} if $|S_n|\le C$ for infinitely many $n$ almost surely. Conversely, we say that the walk is {\em transient} if $\lim |S_n|=\infty$ almost surely.

It was shown in~\cite{bhattacharya2023recurrence} that in the case $a_n=\lfloor n^\alpha\rfloor$ where $\alpha>0$, the walk is always transient.  While the case $\alpha\ge 1$ is somewhat easy due to the increasing size of the steps (see, e.g., \cite{erdos,sarkozi1965erdos,sullivan}), when $0<\alpha<1$ the results are not so trivial, and the case $\alpha\le 1/2$ is especially hard, as the intuition suggests that the behaviour of the walk should depend on $\D(S_n)\asymp n^{\alpha+1/2}$ where $\D(S_n)=\sqrt{\Var(S_n)}$ is the standard deviation, and whether the power of $n$ is smaller than $1$ or not should play an important role. Nonetheless, the walk is transient for both $\alpha<1/2$ and $\alpha\ge 1/2$.

It was left as an open problem to check whether the walk is transient when $a_n=n^\alpha$ (without the integer part), as the methods of~\cite{bhattacharya2023recurrence} to show transience were applicable only to the integer-valued walks. This problem was partially answered in~\cite{BCJ}, where they showed that for $\alpha>1/2$ the walk is, in fact, transient, as conjectured. Some results in higher dimensions were given in~\cite{bhattacharya2025twodimensionalrademacherwalk,englander2025coin}. For an extensive review of the literature on Rademacher random walk and related problems, we also refer the reader to~\cite{BCJ}.

The purpose of this paper is to provide an easy proof of a somewhat simplified\footnote{they have $a_n=n^{\alpha+o(1)}$} version of Theorem~1 in~\cite{BCJ}, using known local limit theorems.

\begin{thm}\label{thm:general}
Let $(a_n)_{n \geq 1}$ be a sequence such that  $a_n  = g(n)$ for $n\ge 1$, where for  $n\ge n_0$ for some $n_0$ the function $g$ is continuously differentiable and satisfies
\begin{align}\label{eqC}
\frac12<p\le\frac{xg'(x)}{g(x)}\le q<\infty,\qquad
g(x)=o(x)
\end{align}
for some $p,q>0$ (and, in particular, $g$ is monotone). Then the associated Rademacher random walk is transient.
\end{thm}
\begin{cor}
Suppose that for large $n$, $a_n=n^\alpha\log^\beta n$ where $\alpha>1/2$, $\beta\in\R$.  Then the associated Rademacher random walk is transient.  
\end{cor}
\begin{proof}
The cases $\alpha> 1$ or $\alpha=1$ and $\beta\ge 0$ follow fairly trivially from the Halász inequality (see e.g.\ Lemma~20 in \cite{BCJ}), The remaining cases where $1/2<\alpha<1$ or $\alpha=1$ and $\beta<0$, follow from Theorem~\ref{thm:general}.   
\end{proof}

We start with some general result. Let $\{S_n\}$ be a sequence of random variables and $\{B_n\}$ be a sequence of positive numbers such that $B_n\to\infty$.
The following statement directly follows from~\cite [Theorem~4.1]{PEL}.
{\begin{lemma}\label{lemmain}
Assume that $S_n/B_n$ converges in distribution to a standard normal random variable. In addition, suppose that for each $D>0$
\begin{equation}\label{cond}
\lim_{T\to\infty}\limsup_{n\to\infty}
\int_{T\le |t|\le DB_n}
\left|\mathbb{E}\exp\left(it\frac{S_n}{B_n}\right)\right|\,dt=0.
\end{equation}
Then for every $C>0$
\begin{equation*}
\lim_{n\to\infty}\left|
B_n\P(S_n\in[-C,C])-\frac{2C}{\sqrt{2\pi}}\right|=0.
\end{equation*}
\end{lemma}}

The following statement is elementary, so we skip its proof.
\begin{lemma}\label{elemma}
For any $x\in\R$, $|\cos(x)|\le e^{-\sin^2(x)/2}$. Also, if
 $|x|\le 1$, then $|\cos(x)|\le e^{-x^2/2}$.
\end{lemma}

Integrating~\eqref{eqC}, we get for some $c_1,c_2>0$
$$
\left(\frac{y}{x}\right)^p\le \frac{g(y)}{g(x)}\le\left(\frac{y}{x}\right)^q
\text{ for $y\ge x\ge n_0$}
\quad \Longrightarrow\quad c_1 n^p\le  g(n) \le c_2 n^q\text{ for all }n.
$$
Let 
$
B_n^2=\Var(S_n)=a_1^2+\dots+a_n^2=\sum_{k=1}^n g(k)^2
$
then, since $g(x)\le g(n)$ for $x\in[n_0,n]$ and for $x\in [n/2,n]$ we have $g(x)\ge 2^{-q} g(n)$,
\[
B_n\asymp\sqrt n\,g(n),\qquad
\sum_n\frac1{B_n}\le \sum_n\frac{c_3}{n^{1/2+p}}<\infty
\]
for some $c_3>0$. Let now
$$
\psi(t)=\left|\E e^{it S_n/B_n}\right|=
\prod_{k=1}^n\left|\cos\left(\frac{t\, a_k}{B_n}\right)\right|.
$$
Our goal is to check condition~\eqref{cond}; due to the symmetry, it suffices to consider the integral over~$t$ only on the interval 
$$
[T,DB_n]=[T,B_n/a_n]\cup [B_n /a_n,D B_n],
$$
On the first part of the above interval as $t a_k/B_n\le t a_n/B_n \le 1$, so by Lemma~\ref{elemma},
\begin{align}\label{eq1}
\psi(t)\le \exp\left(-\frac12\sum_{k=1}^n\frac{t^2\, a_k^2}{B_n^2}\right)=e^{-t^2/2}.    
\end{align}
In the second part of the interval, by the other part of Lemma~\ref{elemma},
\begin{align}\label{eq2}
\psi(t)\le \exp\left(-\frac{M_n(t)}2\right)\qquad\text{where }
M_n(t)=\sum_{k=n/2+1}^n\sin^2\left(\frac{t\, a_k}{B_n}\right)
\end{align}
(we assume that $n$ is even to simplify notations, since the case when $n$ is odd is similar).

\begin{lemma}\label{lem2}
Let $D>0$. Then there exists $\delta>0$ such that for sufficiently large~$n$, we have 
$$  
\inf_{ B_n/a_n\le t\le DB_n}\frac{M_n(t)}{n}\ge \delta.
$$
\end{lemma}
\begin{proof}
Let 
\[
h_n(x)=\frac{g(nx)}{g(n)},\qquad
J_n(u)=\int_{1/2}^1\sin^2\!\bigl(u h_n(x)\bigr)\,dx.
\]
For $1/2\le x\le 1$, since $g$ is monotone and using~\eqref{eqC},
\begin{align}\label{hprime}
0\le h_n'(x)=\frac{(nx)\, g'(nx)}{g(nx)}\cdot \frac{g(nx)}{x\, g(n)}\le q\cdot \frac{g(nx)}{x\, g(n)}\le 2q.    
\end{align}
Therefore, substituting \(y=h_n(x)\), and taking into account that $h_n(1/2)\le 2^{-p}$, $h_n(1)=1$, we get for all $u\ge 1$
\[
J_n(u)\ge \frac1{2q}\int_{2^{-p}}^1\sin^2(uy)\,dy
=\frac1{4q} \int_{2^{-p}}^11-\cos(2uy)\,dy
\ge 2\delta
\]
where \(\delta>0\) is independent of $n$ and $u\ge 1$. Indeed, the last integral is continuous and strictly positive for every \(u\ge1\), and converges to \(1-2^{-p}>0\) as \(u\to\infty\).

Now set $u=\frac{t g(n)}{B_n}$ and $f(x)=\sin^2(uh_n(x))$  Then
\begin{align*}
\left|\frac{M_n(t)}{n}-J_n(u)\right|&=
\left|\frac1n \sum_{k=n/2+1}^n f\left(\frac kn\right)-\int_{1/2}^1 f(x) dx\right|\le 
\sum _{k=n/2+1}^{n}\int_{(k-1)/n}^{k/n}\left|f\left(\frac{k}{n}\right)-f(x)\right|dx
\\ &
\le \sum _{k=n/2+1}^{n}\frac{1}{n}\int _{(k-1)/n}^{k/n}|f'(x)|dx
=\frac{1}{n}\int _{1/2}^{1}|f'(x)|dx
\\ &\le \frac{1}{n}\int_{1/2}^1 u h_n'(x)|\sin (2uh_n(x))|dx
\le \frac1{n}\int_{1/2}^1 2uq dy\le\frac{uq}{n}\le\frac{u}n
\end{align*}
by~\eqref{hprime}.
This gives
\[
\frac{M_n(t)}n \ge J_n(u)-\frac un \ge 2\delta-\frac{Dg(n)}n \ge \delta
\]
for large \(n\), uniformly over $\frac{B_n}{g(n)}\le t\le DB_n.$
\end{proof}

\begin{proof}[Proof of Theorem~\ref{thm:general}]
For sufficiently large $n$ and $T>1$ we have by~\eqref{eq1}, \eqref{eq2} and Lemma~\ref{lem2} 
$$
\left|\mathbb E e^{itS_n/B_n}\right| \le
\begin{cases}
e^{-t^2/2},& |t|\le B_n/g(n),\\[2mm]
e^{-\delta n/2},& B_n/g(n)\le |t|\le DB_n.
\end{cases}
$$
Hence
\begin{align*}
\int_T^{DB_n} \psi(t)dt\le 
\int_T^{\frac{B_n}{g(n)}}\psi(t)dt+
\int_{\frac{B_n}{g(n)}}^{D B_n}\psi(t)dt
&\le \int_T^\infty t e^{-t^2/2}dt+\int_{\frac{B_n}{g(n)}}^{D B_n}e^{-\delta n/2}dt\le e^{-T^2/2}+o(1).
\end{align*}
Thus $
\limsup_{n\to\infty} \int_T^{DB_n} \psi(t)dt \le e^{-T^2/2}
$
yielding the condition~\eqref{cond}. 
Also, since $\max_{k\le n}g(k)/B_n\asymp n^{-1/2}\to 0$, so by the Lindeberg--Feller theorem $S_n/B_n\to\mathcal{N}(0,1)$. Therefore by Lemma~\ref{lemmain}
$$
\P(|S_n|\le C)\le \sqrt{\frac{2}{\pi}}\frac{C+o(1)}{B_n}.
$$
Since $\sum_n 1/B_n<\infty$, the event $\{|S_n|\le C\}$ occurs only finitely often by the Borel-Cantelli lemma. Taking the intersection of this event for every positive integer $C$ finishes the proof,
\end{proof}

\section*{Acknowledgment}
The author used large language models to search the literature, verify minor technical details, and check the proof. The idea behind the proof is entirely the author’s own.

\end{document}